\documentclass[a4paper,11pt,reqno]{amsart}
\usepackage{amsmath} 
\usepackage{amsfonts}
\usepackage{t1enc}
\usepackage{indentfirst}
\usepackage{fancyhdr}
\numberwithin{equation}{section}
\usepackage[margin=2.9cm]{geometry}
\usepackage{hyperref}
\usepackage{enumitem}
\usepackage{etoolbox}
\allowdisplaybreaks[1]

\theoremstyle{plain}
\newtheorem{Th}{Theorem}
\newtheorem{Lemma}[Th]{Lemma}

\newtheorem{Prop}[Th]{Proposition}

\hypersetup{
    colorlinks=true,
    linkcolor=blue,
    citecolor=blue,      
    urlcolor=black,
}

\makeatletter

\patchcmd\subsection
  {.5\linespacing\@plus.7\linespacing}
  {0\linespacing}
  {}
  {}{}

\def\@seccntformat#1{%
  \protect\textup{\protect\@secnumfont
    \ifnum\pdfstrcmp{subsection}{#1}=0 \bfseries\fi
    \csname the#1\endcsname
    \protect\@secnumpunct
  }%
} 

\def\thm@space@setup{\thm@preskip=0pt
\thm@postskip=0pt}

\makeatother

\counterwithout{equation}{section}

\begin{document}

\pagestyle{fancy}
\lhead{}
\chead{}
\rhead{\thepage}
\cfoot{}
\lfoot{}
\rfoot{}
\renewcommand{\headrule}{}

\title{Miscellaneous summation, integration, and transformation formulas} 

\author{Martin Nicholson} 

\begin{abstract}  
This is a discussion of miscellaneous summation, integration and transformation formulas obtained using Fourier analysis. The topics covered are: Series of the form $\sum_{n\in\mathbb{Z}} c_ne^{\pi i \gamma n^2}$; Fusion of integrals, and in particular fusion of $q$-beta integrals related to Gauss-Fourier transform, and a related family of eigenfunctions of the cosine Fourier transform; Summation formulas of the type $\sum_{n\ge 1}\frac{\chi(n)}{n}\,f(n)$ with Dirichlet characters; Trigonometric Fourier series expansion of hypergeometric functions of the argument $\sin^2x$; Modifications of the inverse tangent integral and identities for corresponding infinite products.
\end{abstract}

\clearpage\maketitle
\thispagestyle{empty}
\vspace{-20pt}

\section{Series with general term \texorpdfstring{$c_ne^{\pi i \gamma n^2}$}{cnexp(piign2)}}
 
\subsection{} Consider the series of the form 
\begin{equation}\label{quadratic_series}
    \sum_{n=1}^\infty c_ne^{\pi i n^2/a}.
\end{equation}
Ramanujan has studied \eqref{quadratic_series} with coefficients $c_n$ (\cite{RLN4}, \cite{ramanujan2}, \cite{ramanujan3})
\[
\frac{1}{n^2},\quad \frac{n}{e^{2\pi n}-1}, \quad \frac{1}{\cosh(\pi n/\sqrt{a})},
\]
and many others. A complete description of \eqref{quadratic_series} is given when $a$ is a positive integer by the following fact:
\textit{If $f(x)$ has period $\,1$ and its Fourier series expansion is}
\[
f(x)=\sum_{n\in\mathbb{Z}} c_ne^{2\pi inx},
\]
\textit{then}
\begin{align}\label{even}
    \nonumber\sum_{n\in\mathbb{Z}} c_ne^{\pi i n^2/a}&=\frac{e^{\pi i /4}}{\sqrt{a}}\sum _{r=1}^af\Big(\frac{r}{a}-\frac{1}{2}\Big)e^{-\pi i a(1/2-r/a)^2},\quad a\in\mathbb{N},\\
    &=\frac{e^{\pi i /4}}{\sqrt{a}}\sum _{r=1}^af\Big(\frac{r}{a}\Big)e^{-\pi ir^2/a},\quad \frac{a}{2}\in\mathbb{N}.
\end{align}
\noindent{}This allows one to reduce the series \eqref{quadratic_series} to a finite sum when $a\in\mathbb{N}$ and the function $f(x)$ is known explicitly. This fact is generally well known and the idea behind it is due to Dirichlet \cite{aar}. The case $f(x)=1$ gives the value of the Gauss sum
\begin{equation}\label{gauss_sum}
    \frac{e^{\pi i /4}}{\sqrt{a}}\sum _{r=1}^ae^{-\pi i r^2/a}=1,\quad \frac{a}{2}\in\mathbb{N}.
\end{equation}
As another illustration, we apply it to the sequence $c_0=0$, $c_n=1/2n^{2}$ $(n\neq 0)$, in which case $f(x)$ is known
\[
\pi^2\left(\frac{1}{6}-x+x^2\right)=\sum_{n=1}^\infty \frac{\cos(2\pi nx)}{n^2},\qquad 0\le x\le 1.
\]
Substituting this into \eqref{even} and simplifying the result using \eqref{gauss_sum} yields
\[
\sum_{n=1}^\infty \frac{e^{\pi i n^2/a}}{n^2}=\frac{\pi^2}{6}-\frac{\pi^2}{\sqrt{a}}\sum _{r=1}^a\frac{r}{a}\Big(1-\frac{r}{a}\Big)e^{\pi i/4-\pi ir^2/a},\quad \frac{a}{2}\in\mathbb{N}.
\]
This was recorded in Ramanujan's Lost Notebook (\cite{RLN4}, section 10.3) and first proved in \cite{bct} by different methods, where it was also generalized to sums with $c_n=n^{-2m}$. 
\subsection{} Direct derivation of the fact above based on multisection (also called Simpson dissection, e.g., \cite{aar}) requires calculation of certain Gauss sums (see \cite{bct} for a similar calculation when $c_n=n^{-2m}$). However, a method based on Poisson summation formula avoids calculation of Gauss sums. This method is sketched below. We use Poisson summation formula in the form \cite{aar}
\begin{equation}\label{poisson}
    \sum_{n\in\mathbb{Z}} \frac{g(x+n-0)+g(x+n+0)}{2}=\sum_{n\in\mathbb{Z}}e^{2\pi inx} \int_{-\infty}^{\infty}g(t)e^{-2\pi int}dt.
\end{equation}
Writing the sum as an integral we have, assuming $a\in\mathbb{N}$:
\begin{align}
\frac{e^{-\pi i/4}}{\sqrt{a}}\sum_{n\in\mathbb{Z}} c_ne^{\pi i n^2/a}&=\int_{-\infty}^\infty f(t) e^{-\pi i a t^2}dt\label{integral_rep}\\
\nonumber &=\sum_{n\in\mathbb{Z}}\int_{-1/2}^{1/2} f(t)e^{-\pi i a (t+n)^2}dt\\
\nonumber &=\frac{1}{a}\sum_{n\in\mathbb{Z}} e^{-\pi i a n}\int_{-a/2}^{a/2} f(t/a)e^{-\pi i t^2/a}  e^{-2\pi i n t}dt\\
\nonumber &=\frac{1}{a}\sum _{n=1}^a f\left(\frac{n}{a}-\frac{1}{2}\right)e^{-\pi i a\left(\frac{n}{a}-\frac{1}{2}\right)^2}.
\end{align}
Here, the last equility is the application of \eqref{poisson} to the function
\[
g(t)=\begin{cases}
    f(t/a)e^{-\pi i t^2/a}, &|t|\le a/2,\\
    0, &|t|>a/2.
\end{cases}
\]
The case $a/2\in\mathbb{N}$ can be treated similarly.

\subsection{} The integral representation of the sum in \eqref{integral_rep} is valid for any $a\in\mathbb{R}$. There is the following curious identity related to this integral representation:
\textit{If the real numbers $\alpha$ and $\beta$ are such that $\alpha\beta \in\mathbb{Z}$ and the functions $f(x)$ and $g(x)$ have period $1$, then}
\[
\int_{-\infty}^\infty f(\alpha x)g(\beta x) e^{\pi i x^2}dx=e^{-\pi i/4}\int_{-\infty}^\infty f(\alpha x) e^{\pi i x^2}dx \int_{-\infty}^\infty g(\beta x) e^{\pi i x^2}dx.
\]
\noindent{}This is due to the fact that if one writes
\[
f(\alpha x)g(\beta x)=\sum_{r\in\mathbb{Z}} a_re^{2\pi i\alpha rx}\sum_{s\in\mathbb{Z}} b_se^{2\pi i\beta sx},
\]
and calculates the resulting integral over $x$ as
\[
\int_{-\infty}^\infty e^{\pi i x^2+2\pi ix(\alpha r+\beta s)}dx= e^{-\pi i/4-\pi i\alpha^2 r^2-\pi i\beta^2 s^2}\cdot e^{-2\pi i\alpha\beta rs},
\]
then the term mixing the sums over $r$ and $s$ is
\[
e^{-2\pi i\alpha\beta rs}\equiv 1
\]
when $\alpha\beta$ is an integer, so the sums over $r$ and $s$ decouple. A less trivial type of factorization is considered in section \ref{q_integrals}, where we discuss the fusion of $q$-beta integrals.

\subsection{}  An entry in Ramanujan's Lost Notebook (\cite{RLN2}, Entry 1.7.16) states that: 
\textit{Let $a$ and $b$ be any complex numbers, and suppose that $|xy|<1$. If}
\[
\phi(a,x,y)=\sum_{n=0}^\infty\frac{a^nx^{n(n+1)/2}}{(xy;xy)_n},
\]
\textit{then}
\[
\phi(a,x,y)\phi(b,y,x)=\sum_{n=0}^\infty\frac{(ax+by^n)(ax^2+by^{n-1})\ldots (ax^n+by)}{(xy;xy)_n}.
\]
In the limit $x,y\to 1^{-}$, this identity reduces to
\begin{equation}\label{1716}
    \left|\sum_{n=0}^\infty\frac{e^{i\gamma n^2}}{n!}\right|^2=\sum_{n=0}^\infty\frac{2^n}{n!}(\cos\gamma n)^n,\quad \gamma\in\mathbb{R},
\end{equation}
which can be easily verified directly using binomial theorem.

\subsection{} Another formula similar to \eqref{1716} is 
\begin{align*}
&\left|\sum_{n\in\mathbb{Z}}\frac{e^{\pi  i  \gamma  n^2}}{p^2+4\pi ^2 n^2}\right|^2=\frac{1}{p \sinh \frac{p}{2}}+8p\sinh\tfrac{p}{2}\\
&+\sum_{n\in\mathbb{Z}}\frac{2}{p^2+\pi^2 n^2}\left[\cos \left(\pi  \gamma n^2\right) \cosh \left(p \left\{\gamma n\right\}-\frac{p}{2}\right)-\frac{p}{\pi n}{\sin\left(\pi\gamma n^2\right) \sinh\left(p \left\{\gamma  n\right\}-\frac{p}{2}\right)}\right],
\end{align*}
where $\gamma\in\mathbb{R}$, and $\left\{x\right\}$ denotes the fractional part of the real number $x$. Its continuous counterpart is
\[
\left|\int_0^\infty\frac{e^{i\alpha x^2}}{x^2+1}dx\right|^2=\pi\int_0^\infty\frac{e^{-2\alpha x}}{x^2+4}\left(\cos\alpha x^2+\frac2x \sin\alpha x^2\right)dx.
\]
These type of integrals are related to Owen's T-function \cite{owen}
\[
T(h,a)=\frac{1}{2\pi}\int_0^a\frac{e^{-h^2(1+x^2)/2}}{1+x^2}\, dx.
\]
Closely related is also entry $7.4.24$ in \cite{as}.

\subsection{} If $\alpha\beta=1$, $q=e^{-\pi\alpha}$, then
\begin{align}\label{elliptic}
    \nonumber\frac{1}{2\alpha}
&\left|\sum_{n\in\mathbb{Z}}\frac{e^{\pi i \gamma n^2}}{\cosh(\pi\beta n)}\right|^2-\sum_{n\in\mathbb{Z}}\frac{\alpha n}{\sinh(\pi \alpha n)}\\&=2\frac{(q^2;q^2)_\infty^4}{(q;q)_\infty^2}\sum_{n=1}^\infty\frac{(q^{1+2n\gamma},q^{1-2n\gamma};q^2)_\infty}{(q^{2+2n\gamma},q^{2-2n\gamma};q^2)_\infty}\,\frac{\sin(\pi \gamma n^2)}{\sinh(\pi \alpha\gamma n)\sinh(\pi \beta n)}.
\end{align}
Continuous analog of \eqref{elliptic} is a relation between Mordell type integrals
\[
\int_0^\infty\frac{\sin(\gamma x^2)}{\sinh(\pi x)\sinh(\gamma x)}dx=\left|\int_0^\infty\frac{e^{i\gamma x^2}}{\cosh(\pi x)}dx\right|^2.
\]
In terms of elliptic integrals $K(k)=\int_0^{\pi/2}\frac{d\theta}{\sqrt{1-k^2\sin^2\theta}}$, $E(k)=\int_0^{\pi/2}{\sqrt{1-k^2\sin^2\theta}}\,d\theta$, $k'=\sqrt{1-k^2}$, $K'=K(k')$, $E'=E(k')$, formula \eqref{elliptic} takes the form
\[
\frac{4E'}{\pi}+\sum_{n=1}^\infty\frac{8i}{\operatorname{sn}(2i\gamma K'n,k)}\frac{\sin(\pi \gamma n^2)}{\sinh(\pi K n/K')}=\frac{\pi}{K'}
\left|\sum_{n\in\mathbb{Z}}\frac{e^{\pi i \gamma n^2}}{\cosh(\pi K n/K')}\right|^2.
\]
Series with elliptic functions are studied in the theory of elliptic hypergeometric functions \cite{gr2}. The sum
\[
f(\alpha)=\sum_{n\in\mathbb{Z}}\frac{e^{\pi i \gamma n^2}}{\cosh(\pi\alpha n)}
\]
with $\gamma=\alpha^2$ was studied in \cite{ramanujan2}, where it was shown that certain linear combination of $f(\alpha)$ and complex conjugate of $f(1/\alpha)$ is a modular function of the variable $i\alpha$ .

\section{Fourier-Gauss transform and fusion of integrals}

\subsection{} Two integrals (\cite{RLN1}, Entries 16.1.4-5) 
\begin{equation}\label{beta1}
    \int_{-\infty}^\infty e^{-x^2+2imx}(a\sqrt{q}\,e^{2k x},b\sqrt{q}\,e^{-2k x};q)_\infty\, dx=\sqrt{\pi}\,e^{-m^2}\frac{(ab;q)_\infty}{(-ae^{2ikm},-be^{-2ikm};q)_\infty},\quad q=e^{-2k^2},
\end{equation}
\begin{equation}\label{beta2}
    \int_{-\infty}^\infty e^{-x^2+2imx}\frac{1}{(-ae^{2ik x},-be^{-2ik x};q)_\infty}\,dx=\sqrt{\pi}\,e^{-m^2}\frac{(a\sqrt{q}\,e^{-2k m},b\sqrt{q}\,e^{2k m};q)_\infty}{(ab;q)_\infty},\quad q=e^{-2k^2},
\end{equation}
involving infinite $q$-products, were proved by Askey using $q$-binomial theorem \cite{askey}. Pastro gave an alternative proof \cite{pastro}. Stokman fused these integrals into a single integral with $4$ parameters \cite{stokman}. Note that \eqref{beta1} and \eqref{beta2} are not independent: If one of them is given, then the other follows by inverse Fourier transform.

Let us define Fourier-Gauss transform of a function $f(x)$ by the formula
\[
\int_{-\infty}^\infty f(x)\,e^{-x^2+2imx}\,dx.
\]
It is a special case of the so called Fourier-Bros-Iagolnitzer transform \cite{folland}.
Using the value of the gaussian integral (\cite{gr1}, Entries 4.133.1-2)
\begin{equation}\label{gauss_integral}
    \int_{-\infty}^\infty e^{-\frac{x^2}{4\gamma}+\beta x+iax}dx=2\sqrt{\pi\gamma}\, e^{\gamma(\beta^2-a^2)+2i\gamma a\beta}.
\end{equation}
we find that the action of Fourier-Gauss transform on the exponential function is given by
\begin{equation}\label{exp1}
    \int_{-\infty}^\infty e^{2knx}\,e^{-x^2+2imx}\,dx=\sqrt{\pi}\,e^{-m^2}\cdot e^{k^2n^2+2iknx},
\end{equation}
\begin{equation}\label{exp2}
    \int_{-\infty}^\infty e^{2iknx}\,e^{-x^2+2imx}\,dx=\sqrt{\pi}\,e^{-m^2}\cdot e^{-k^2n^2-2knx}.
\end{equation}
Now, let $u_n(k)$ be an absolutely summable sequence, (i.e., a sequence such that $\sum_{n\in\mathbb{Z}} |u_n(k)|<\infty$), and define the functions $\phi_u(x,k)$, $\psi_u(x,k)$ related to this sequence according to
\begin{equation}\label{phi}
    \phi_u(x,k)=\sum_{n\in\mathbb{Z}} u_n(k)e^{-k^2n^2+2k n x},
\end{equation}
\begin{equation}\label{psi}
    \psi_u(x,k)=\sum_{n\in\mathbb{Z}} u_n(k)e^{2ik n x}.
\end{equation}
It follows from \eqref{exp1} and \eqref{exp2}, that the functions \eqref{phi}--\eqref{psi} are related to each other by Fourier-Gauss transform:
\[
\int_{-\infty}^\infty e^{-x^2+2imx}\phi_u(x,k)\, dx=\sqrt{\pi}\,e^{-m^2}\psi_u(m,k),
\]
\[
\int_{-\infty}^\infty e^{-x^2+2imx}\psi_u(x,k)\, dx=\sqrt{\pi}\,e^{-m^2}\phi_u(-m,k).
\]
This property of Fourier-Gauss transform was used in \cite{atakishiyev} to relate different families of orthogonal polynomials to each other.

Thus, \eqref{beta1} and \eqref{beta2} give an example of pair of complementary functions \eqref{phi}--\eqref{psi} for a particular choice of the sequence $u_n(k)$:
\begin{equation}\label{phi2}
    \phi_u(x,k)=(a\sqrt{q}\,e^{2k x},b\sqrt{q}\,e^{-2k x};q)_\infty,
\end{equation}
\begin{equation}\label{psi2}
    \psi_u(x,k)=\frac{(ab;q)_\infty}{(-ae^{2ikm},-be^{-2ikm};q)_\infty}.
\end{equation}
For a pair of complimentary functions containing 3 parameters, see Entry 1.7.19 in \cite{RLN2}.

\subsection{}\label{q_integrals} Another property of the pair of functions \eqref{phi}--\eqref{psi} is given below. 
\begin{Prop}\label{fusion} Let $k\in\mathbb{R}$. Let $u_n(k)$, $v_n(k)$ be absolutely summable sequences. For $x\in\mathbb{R}$ and $u_n(k)$ define $\phi_u(x,k)$, $\psi_u(x,k)$ by \eqref{phi}, \eqref{psi} and similarly $\phi_v(x,k)$, $\psi_v(x,k)$ for $v_n(k)$. Then for $m\in\mathbb{R}$
\[
\int_{-\infty}^\infty e^{-x^2+2imx}\phi_u(x,k)\psi_v(x,\pi/k)~dx=\sqrt{\pi}\,e^{-m^2}\phi_v(-m,\pi/k)\psi_u(m,k).
\]
\end{Prop}

\noindent{}{\it{Proof.}} We write
\[
\phi_u(x,k)\psi_v(x,\tilde{k})=\sum_{r\in\mathbb{Z}} u_r(k)e^{-k^2r^2+2k r x}\sum_{s\in\mathbb{Z}} v_s(\tilde{k})e^{2i\tilde{k} s x}.
\]
Then by \eqref{gauss_integral}
\[
\int_{-\infty}^\infty e^{-x^2+2imx+2krx+2i\tilde{k}sx}dx=\sqrt{\pi}\, e^{k^2r^2-(m+\tilde{k}s)^2+2ikmr+2ik\tilde{k}rs}.
\]
Thus, after integration, the mixing of the sums over $r$ and $s$ occurs through the factor $e^{2ik\tilde{k}rs}$ only. This factor equals $1$ when $k\tilde{k}=\pi$ (more generally, when $k\tilde{k}/\pi\in\mathbb{N}$), in which case the sums over $r$ and $s$ decouple. The sum over $r$ results in the factor $\phi_v(-m,\pi/k)$, while the sum over $s$ results in the factor $\psi_u(m,k)$. \qed

Coupling Proposition \ref{fusion} with equations \eqref{phi2} and \eqref{psi2} one obtains that for $0\le |a|,|b|,|c|,|d|<1$
\begin{equation}\label{fused_q_integral}
    \int_{-\infty}^\infty e^{-x^2}\frac{(a\sqrt{q}\,e^{-2k x},b\sqrt{q}\,e^{2k x};q)_\infty}{(-ce^{-2\pi i x/k},-de^{2\pi ix/k};Q)_\infty}\, dx=\sqrt{\pi}\,\frac{(ab;q)_\infty(c\sqrt{Q},d\sqrt{Q};Q)_\infty}{(-a,-b;q)_\infty(cd;Q)_\infty},
\end{equation}
where $Q=e^{-2\pi^2/k^2}$. This is Proposition $3.1$ in \cite{stokman}. Integral \eqref{fused_q_integral} is fusion of integrals \eqref{beta1} and \eqref{beta2}. Thus Proposition \ref{fusion} is the fusion property formulated in terms of Fourier-Gauss transform. 

\subsection{} \textit{A family of eigenfunctions of the cosine Fourier transform.} Define
cosine Fourier transform of the function $f(x)$ by the formula
\[
F_c(y)=\sqrt{\frac{2}{\pi}}\int_0^\infty f(x)\cos (xy) \,dx.
\]
It follows from Proposition \ref{fusion}, that for a suitable set of parameters it is possible to obtain functions that are their own Fourier transforms (eigenfuntions of the Fourier transform): \textit{If}
\[
f(x)=e^{-x^2/2}\sum_{r=0}^\infty u_r\,e^{-\pi r^2}\cosh\big(\sqrt{2\pi}\,rx\big)\sum_{s=0}^\infty u_s\cos\big(\sqrt{2\pi}\,sx\big),
\]
\[
g(x)=e^{-x^2/2}\sum_{r=0}^\infty u_r\,e^{-\pi r^2}\sinh\big(\sqrt{2\pi}\,rx\big)\sum_{s=0}^\infty u_s\sin\big(\sqrt{2\pi}\,sx\big),
\]
\textit{then}
\[
F_c(x)=f(x),\quad G_c(x)=-g(x).
\]
For example
\[
f(x)=e^{-x^2/2}\frac{(a\sqrt{q}\,e^{-\sqrt{2\pi} x},a\sqrt{q}\,e^{\sqrt{2\pi} x};q)_\infty}{(-ae^{-\sqrt{2\pi}ix},-ae^{\sqrt{2\pi}ix};q)_\infty},\quad q=e^{-2\pi},
\]
satisfies $F_c(x)=f(x)$.

\subsection{} According to generalized Plancherel's theorem, it follows from Proposition \ref{fusion} that
\begin{align*}
    \int_{-\infty}^\infty e^{-2x^2}\phi_u(x,k)\psi_v(x,\pi/k)&\phi_\xi(x,k)\psi_\eta(x,\pi/k)\,dx\\
    &=\int_{-\infty}^\infty e^{-2x^2}\phi_v(-x,\pi/k)\psi_u(x,k)\phi_\eta(x,\pi/k)\psi_\xi(-x,k)\,dx.
\end{align*}
Calculation of such integrals when $\phi$ is a function of the type \eqref{phi}  was discussed in \cite{suslov}, $r=2$ case of equation (4.1). Integrals in the right hand side of equation (4.1) in the article \cite{suslov} are calculated using a trick due to Bailey \cite{bailey}. The result is an integral with 8 parameters expressed as a sum of two products of two basic bilateral series ${}_2\psi_2$ with bases ${q}=e^{-2k^2}$ and $Q=e^{-2\pi^2/k^2}$, respectively. In fact, the $r=2$ case of equation (4.1) is a generalization of modular transformation formula for Appell-Lerch sums \cite{gm}. This is because Appell-Lerch sums have representation in terms of ${}_2\psi_2$ \cite{choi}.

A particular case of the said general 8 parameter integral is: \textit{Let $k>0$, $q=e^{-2k^2}$, $Q=e^{-2\pi^2/k^2}$, $a=e^{2k\alpha}$, $b=e^{2k\beta}$, $\operatorname{arg}\,(ab)<\pi$. Then}
\begin{align}\label{al}
\nonumber\int_{-\infty}^{\infty} &e^{-x^2}\frac{\left(q a e^{2k x},q b e^{2kx};q\right)_{\infty }}{\left(-\sqrt{q}ab e^{2k x},-\sqrt{q}e^{-2k x}/ab;q\right)_\infty} \, dx=\frac{\sqrt{\pi }}{\left(-\sqrt{q}a,-\sqrt{q}/a,-\sqrt{q}b,-\sqrt{q}/b;q\right)_\infty}\\&-\frac{\pi}{k^2}\sqrt{\frac{\pi a}{q b}}\frac{e^{i \pi  (\alpha +\beta )/k-2 \alpha  \beta-\beta^2}}{\left(q^3Q\right)^{1/8}(q;q)_\infty^4}
\sum_{n\in\mathbb{Z}} \frac{b^n q^{{n^2}/{2}}}{1+a q^{n-{1}/{2}}}
\sum_{n\in\mathbb{Z}}\frac{e^{-{2 \pi  i\alpha n}/{k}}{Q}^{n^2/2}}{1+e^{{2 \pi  i\beta}/{k}}{Q}^{n-1/2}}.
\end{align}
Using properties of Appell-Lerch sums one can deduce from \eqref{al} that the integral
\[
I(\alpha,\beta,k)=q^{1/12}e^{(\alpha +\beta )^2/2}\int_{-\infty}^{\infty}\frac{\big(q e^{2k (x+\alpha)},q e^{2k(x+\beta)};q\big)_{\infty }}{\big(-\sqrt{q}e^{2k(x+\alpha+\beta)},-\sqrt{q}e^{-2k (x+\alpha+\beta)};q\big)_\infty} \, e^{-x^2} dx,\qquad q=e^{-2k^2},
\]
has the following non-trivial symmetries
\[
I(\alpha,\beta,k)=I(-i\alpha,i\beta,\pi/k)=I(-\alpha,-\beta,k).
\]

\section{Series with general term \texorpdfstring{$\frac{\chi_k(n)}{n}f(\alpha n)$}{chi(n)f(an)/n}}\label{sums_with_dirichlet_characters}

\subsection{} The formal identity
\begin{equation}\label{formal1}
    \sum_{n=0}^\infty\frac{(-1)^{n}\left\{\varphi(2n+1)+\varphi(-2n-1)\right\}}{2n+1}=\frac{\pi}{2}\varphi(0),
\end{equation}
was considered by Ramanujan. It was put on a rigorous ground in \cite{RNI}, p.95-96. One of the conditions for the validity of formulas of this type is that $\varphi(z)$ must be an entire function with bounded growth rate, i.e. $|\varphi(z)|<Ce^{{\pi |z|}/{2}}$, $|z|\to\infty$. Berndt proves this theorem by applying contour integration to the function
\[
\frac{f(z)}{z\cos(\pi z/2)}.
\]
One may consider
\[
\frac{f(z)}{z\sin(\pi (z-s))\sin(\pi (z+s))},
\]
where $f(z)$ is an even entire function with bounded growth rate $|f(z)|<Ce^{2\pi |z|}$, $|z|\to\infty$, and by a slight change of Berndt's argument deduce
\begin{equation}\label{formal2}
    \sum_{n\in\mathbb{Z}}\frac{f\left(n+s\right)}{n+s}=\pi \cot(\pi s) f(0).
\end{equation}
When $\varphi$ is an even function, \eqref{formal1} becomes
\begin{equation}\label{formal1even}
    \sum_{n=0}^\infty\frac{(-1)^{n}\varphi(2n+1)}{2n+1}=\frac{\pi}{4}\varphi(0).
\end{equation}
If we set $s=1/4$ in \eqref{formal2} and rearrange the terms of the series (assuming it is justified) we obtain
\begin{equation}\label{formal4}
    \sum_{n=0}^\infty\left(\frac{f(n+1/4)}{4n+1}-\frac{f(n+3/4)}{4n+3}\right)=\frac{\pi}{4}f(0).
\end{equation}
This is the same as \eqref{formal1even} under the identification $f(x)=\varphi(4x)$. An interesting discussion of \eqref{formal1} is contained in \cite{gosper}, where it is generalized to a sum over roots of Bessel functions.

\subsection{} According to Paley-Wiener theorem, the set of entire functions $f(z)$ with bounded growth rate $|f(z)|<Ce^{\sigma |z|}$, $|z|\to\infty$, that are square integrable on the real axes, coincides with the set of functions $B_\sigma$ that admit representation
\[
f(x)=\int_{-\sigma}^{\sigma} e^{ixt}\varphi(t)dt,
\]
with $\varphi(t)$ square integrable on $[-\sigma,\sigma]$. In other words, $f(x)$ is a function whose Fourier spectrum is limited to the band $[-\sigma,\sigma]$.
Thus \eqref{formal2} is also valid for even band limited functions $f(x)\in B_{2\pi}$. However, it is instructive to give a proof of \eqref{formal2} for band-limited functions which is independent of the Paley-Wiener theorem. Our starting point is the fact that the series 
\begin{equation}\label{piecewise_continuous_function}
    S(x)=\sum_{n\in\mathbb{Z}}\frac{\cos\left(2\pi x(n+s)\right)}{n+s}=\pi\cot(\pi s),\qquad x,s\in (0,1),
\end{equation}
represents a piecewise continuous function (\cite{gr1}, entries 1.445.5-6; also, section \ref{fourier_series}). 

Let $f(x)$ be an even band-limited function
\[
F_c(y)=0, \quad |y|>2\pi.
\]
Multiplying \eqref{piecewise_continuous_function} by $\sqrt{8\pi}F_c(2\pi x)$ and integrating with respect to $x$ from $0$ to $1$, using
\[
\sqrt{8\pi}\int_0^1 F_c(2\pi x)\cos\left(2\pi x(n+s)\right)dx=\sqrt{\frac{2}{\pi}}\int_0^\infty F_c(y)\cos\left(y(n+s)\right)dy=f(n+s),
\]
\[
\sqrt{8\pi}\int_0^1 F_c(2\pi x)dx=\sqrt{\frac{2}{\pi}}\int_0^\infty F_c(y)dy=f(0),
\]
we recover \eqref{formal2}. Formula \eqref{formal1} for band limited functions was studied in \cite{zayed}, equation (2.5).

\subsection{}\label{band_limited} General sampling theorems twisted by Dirichlet characters were studied in \cite{klusch}. Note that \eqref{formal4} may be written as
\[
\sum_{n=1}^\infty\frac{\chi_4(n)}{n}f(n/4)=\frac{\pi}{4}f(0),
\]
where $\chi_4(n)=\sin\frac{\pi n}{2}$ is the odd Dirichlet character modulo $4$.
Notice that the coefficient $\pi/4$ on the left hand side is the value of the Dirichlet $L$-series for the character $\chi_4(n)$. We wish to generalize this formula for other odd Dirichlet characters.

Let $\chi_k(n)$ be any $k$-periodic odd sequence
\begin{equation}\label{odd_periodic_sequence}
    \chi_k(n)=\chi_k(n+k),\quad \chi_k(-n)=-\chi_k(n).
\end{equation}
In particular, equations \eqref{odd_periodic_sequence} are satisfied by odd Dirichlet character $\chi_k(n)$ modulo $k$. Let $L_k(s,\chi_k)$ be the corresponding Dirichlet $L$-series
\[
L_k(s,\chi_k)=\sum_{n=1}^\infty\frac{\chi_k(n)}{n^s}.
\]

Now, we set $s=l/k$ in \eqref{formal2}, where $l$ is an integer such that $0<l<k$, and then rewrite the right hand side using partial fractions expansion of the cotangent function
\[
\sum_{n=0}^\infty \left(\frac{f((nk+l)/k)}{nk+l}-\frac{f((nk+k-l)/k)}{nk+k-l}\right)=\sum_{n=0}^\infty \left(\frac{1}{nk+l}-\frac{1}{nk+k-l}\right)f(0).
\]
It immediately follows from this that
\begin{Prop}
  Let $f(z)\in B_{2\pi}$ be an even band limited function. Let the sequence $\chi_k(n)$ be a $k$-periodic odd sequence \eqref{odd_periodic_sequence}. Then  
\begin{equation}\label{general_sum}
    \sum_{n=1}^\infty\frac{\chi_k(n)}{n}f(n/k)=\sum_{n=1}^\infty\frac{\chi_k(n)}{n}f(0).
\end{equation}
In particular, if $\chi_k(n)$ is an odd Dirichlet character mod $k$, then
\begin{equation}\label{dirichlet_general}
    \sum_{n=1}^\infty\frac{\chi_k(n)}{n}f(an/k)=L_k(1,\chi_k)f(0),\qquad 0\le a\le 1.
\end{equation}
\end{Prop}

As an example, consider the function
\[
f(x)=\frac{\sin\sqrt{b^2+(2\pi x)^2}}{\sqrt{b^2+(2\pi x)^2}}.
\]
It follows from Entry 3.876.1 in \cite{gr1}
\[
\int_0^\infty\frac{\sin\sqrt{b^2+x^2}}{\sqrt{b^2+x^2}}\cos(xy)\, dx=\begin{cases}
    \frac{\pi}{2}J_{0}\left(b\sqrt{1-y^2}\right),& 0<y<1,\\
    0,& y>1,
\end{cases}
\]
that $f\big(\tfrac{x}{2\pi}\big)\in B_1$, and thus $f(x)\in B_{2\pi}$. Hence, it follows from \eqref{dirichlet_general} that
\begin{equation}\label{dirichlet3}
    \sum_{n=1}^\infty\frac{\chi_k(n)}{n}\frac{\sin\sqrt{b^2+(a n/k)^2}}{\sqrt{b^2+(an/k)^2}}=L_k(1,\chi_k)\frac{\sin b}{b},\qquad 0\le a\le 2\pi.
\end{equation}
The $k=4$, $a=2\pi$ case of \eqref{dirichlet3} was recorded in \cite{gosper} as
\[
\sum_{n=0}^\infty\frac{(-1)^n}{2n+1}\frac{\sin\sqrt{b^2+\pi^2 (n+1/2)^2}}{\sqrt{b^2+\pi^2 (n+1/2)^2}}=\frac{\pi}{2}\frac{\sin b}{b}.
\]

\subsection{} We continue to study the sum \eqref{formal1even}. Define the sine Fourier transform $F_s$ of the function $f$ as
\[
F_s(y)=\sqrt{\frac{2}{\pi}}\int_0^\infty f(x)\sin (xy) dx.
\]
Note the Poisson summation formula for the odd Dirichlet character $\chi_4(n)=\sin\tfrac{\pi n}{2}$ modulo $4$ \cite{titchmarsh}
\[
\sqrt{\alpha}\sum_{n=1}^\infty\chi_4(n)f(\alpha n)=\sqrt{\beta}\sum_{n=1}^\infty\chi_4(n)F_s(\beta n),\quad \alpha\beta=\frac{\pi}{2}.
\]

\begin{Prop}\label{summation_formula1} If $\alpha\beta=\frac{\pi}{2}$, then
\[
\sum_{n=1}^\infty\frac{\chi_4(n)}{n}f\left(\alpha n\right)=\beta\sqrt{\frac{\pi}{2}}\int_0^\infty (-1)^{\lfloor t+\frac12\rfloor}F_c(2\beta t)dt.
\]
\end{Prop}
\noindent{}{\it{Proof.}} Our proof follows the proof of formula (1.6) in the article \cite{gosper}. We start from the Fourier series expansion
\begin{equation}\label{fourier_series_chi4}
    (-1)^{\lfloor t+\frac12\rfloor}=\frac{4}{\pi}\sum_{n=1}^\infty\frac{\chi_4(n)}{n}\cos(\pi nt).
\end{equation}
Multiplying it with $F_c(2\beta t)$ and integrating we find
\begin{align*}
    \beta\sqrt{\frac{\pi}{2}}\int_0^\infty (-1)^{\lfloor t+\frac12\rfloor}F_c(2\beta t)dt&=2\beta\sum_{n=1}^\infty\frac{\chi_4(n)}{n}\sqrt{\frac{2}{\pi}}\int_0^\infty F_c(2\beta t)\cos(\pi nt)\,dt\\
    &=\sum_{n=1}^\infty\frac{\chi_4(n)}{n}f\left(\frac{\pi n}{2\beta}\right),
\end{align*}
as required. \qed

\subsection{} Note the Poisson summation formula for the odd Dirichlet character $\chi_3(n)=\left(\frac{n}{3}\right)$ modulo $3$ \cite{titchmarsh}
\[
\sqrt{\alpha}\sum_{n=1}^\infty\chi_3(n)f(\alpha n)=\sqrt{\beta}\sum_{n=1}^\infty\chi_3(n)F_s(\beta n),\qquad \alpha\beta=\frac{2\pi}{3}.
\]

\begin{Prop}\label{summation_formula2} Let $\alpha$ and $\beta$ be real numbers such that $\alpha\beta=\tfrac{2\pi}{3}$, and let $\chi_3(n)$ denote the odd Dirichlet character mod $3$. Then
\[
\sum_{n=1}^\infty\frac{\chi_3(n)}{n}f\left(\alpha n\right)=\beta\sqrt{\frac{2\pi}{3}}\int_0^\infty \epsilon(t) F_c(3\beta t)dt,
\]
where
\[
\epsilon(t)=\left\{\begin{array}{cc}
    -2, &  \text{when } t\in\left(n+\tfrac{1}{3},n+\tfrac{2}{3}\right), n\in \mathbb{Z}\\
    1, & \text{otherwise}
\end{array}\right\}.
\]
\end{Prop}

\noindent{}{\it{Proof.}} Analogous to the proof of Proposition \ref{summation_formula1}, using the Fourier series
\[
\epsilon(t)=\frac{3\sqrt{3}}{\pi}\sum_{n=1}^\infty\frac{\chi_3(n)}{n}\cos(2\pi nt)
\]
instead of \eqref{fourier_series_chi4}.\qed

Propositions \ref{summation_formula1} and  \ref{summation_formula2} can be generalized to other Dirichlet characters. However, the functions $\epsilon(t)$ will become more convoluted.

\subsection{}\label{examples}
We consider several examples illustrating Propositions \ref{summation_formula1} and \ref{summation_formula2}.

\noindent{}{\it{Example 1.}} The function
\begin{equation}\label{self_reciprocal_function1}
    h(x)=\frac{1}{\cosh\left(\sqrt{\frac{\pi}{2}}x\right)}
\end{equation}
is self-reciprocal, i.e. it satisfies $H_c(x)=h(x)$ \cite{titchmarsh}. We have (\cite{RNII}, chapter 14, Entry 15)
\begin{equation}\label{bvp6}
    \sum_{n=1}^\infty\frac{\chi_4(n)}{n}h(\alpha n)+\sum_{n=1}^\infty\frac{\chi_4(n)}{n}h(\beta n)=\frac{\pi}{4},\qquad \alpha\beta=\frac{\pi}{2}.
\end{equation}
Combining the above facts with Proposition \ref{summation_formula1} we arrive at
\begin{equation}\label{example1}
    \alpha\int_0^\infty \frac{(-1)^{\lfloor t+\frac12\rfloor}}{\cosh(\pi \alpha t)}\,dt+\beta\int_0^\infty \frac{(-1)^{\lfloor t+\frac12\rfloor}}{\cosh(\pi \beta t)}\,dt=\frac12,\qquad \alpha\beta=1.
\end{equation}
This formula is in fact an alternative form of the Jacobi's imaginary transformation for the modular angle. When $\alpha=\beta=1$ we get the closed form
\[
\int_0^\infty \frac{(-1)^{\lfloor t+\frac12\rfloor}}{\cosh(\pi t)}\,dt=\frac14.
\]

\noindent{}{\it{Example 2.}} Consider the self-reciprocal function \cite{titchmarsh}
\[
h(x)=\frac{1}{1+2\cosh\left(\sqrt{\frac{2\pi}{3}}\,x\right)}.
\]
It is known that (see footnote on page 4 in \cite{grinberg})
\[
\sum_{n=1}^\infty\frac{(-1)^{n-1}}{n}\frac{\sinh\frac{\pi ny}{a}}{\sinh\frac{\pi nb}{a}}\sin\frac{\pi n x}{a}+\sum_{n=1}^\infty\frac{(-1)^{n-1}}{n}\frac{\sinh\frac{\pi nx}{b}}{\sinh\frac{\pi na}{b}}\sin\frac{\pi n y}{b}=\frac{\pi xy}{ab}.
\]
Simplifying further by setting $x=a/3$ and $y=b/3$ one obtains
\begin{equation}\label{bvp3}
    \sum_{n=1}^\infty\frac{1}{n}\left(\frac{n}{3}\right)\frac{1}{1+2\cosh\frac{2\pi bn}{3a}}+\sum_{n=1}^\infty\frac{1}{n}\left(\frac{n}{3}\right)\frac{1}{1+2\cosh\frac{2\pi an}{3b}}=\frac{\pi}{9\sqrt{3}}.
\end{equation}
Hence
\[
\sum_{n=1}^\infty\frac{\chi_3(n)}{n}h(\alpha n)+\sum_{n=1}^\infty\frac{\chi_3(n)}{n}h(\beta n)=\frac{\pi}{9\sqrt{3}},\qquad \alpha\beta=\frac{2\pi}{3}.
\]
Thus, it follows from Proposition \ref{summation_formula2} that
\[
\alpha\int_0^\infty \frac{\epsilon(t)\,dt}{1+2\cosh(2\pi \alpha t)}+\beta\int_0^\infty \frac{\epsilon(t)\,dt}{1+2\cosh(2\pi\beta t)}=\frac{1}{6\sqrt{3}},\qquad \alpha\beta=1.
\]
The case $\alpha=\beta=1$ results in the closed form
\[
\int_0^\infty \frac{\epsilon(t)\,dt}{1+2\cosh(2\pi t)}=\frac{1}{12\sqrt{3}}.
\]

\noindent{}{\it{Example 3.}} Another integral evaluation is
\begin{equation}\label{example4}
    \int_0^\infty \frac{\arctan \left(e^{-\pi \alpha t}\right)-\arctan \left(e^{-\pi\beta t}\right)}{t}\,(-1)^{\lfloor t+\frac12\rfloor}\,dt=\frac{\pi}{8}\log \frac{\beta}{\alpha},\quad \alpha\beta=1.
\end{equation}
This formula is more interesting than the integrals considered in examples 1 and 2 because in each region of continuity of the function $(-1)^{\lfloor t+\frac12\rfloor}$, the integrand does not have a closed form anti-derivative. Differentiation with respect to $\alpha$ under the integral sign reduces \eqref{example4} to \eqref{example1}. A similar formula may be deduced from Example 2. Note also that according to Frullani's formula
\[
\int_0^\infty \frac{\arctan \left(e^{-\pi \alpha t}\right)-\arctan \left(e^{-\pi\beta t}\right)}{t}\,dt=\frac{\pi}{4}\log \frac{\beta}{\alpha}.
\]

\subsection{} Though not a direct consequence of Propositions \ref{summation_formula1} and \ref{summation_formula2}, the following integral evaluation is based on the same principles
\begin{equation}\label{integral_of_piecewise_function}
    \int_0^\infty \frac{(-1)^{\lfloor t+\frac12\rfloor}}{\cosh(\pi t)}\,e^{\pi i t^2}dt=\frac{1}{2\sqrt{2}}.
\end{equation}
Using the formula (\cite{gr1}, entries 3.989.3-4)
\[
\int_0^\infty \frac{\cos(\pi a t)}{\cosh(\pi t)}\,e^{\pi i t^2}dt=\frac{ie^{-\pi ia^2/4}+e^{-\pi i/4}}{2\cosh\frac{\pi a}{2}},
\]
\eqref{integral_of_piecewise_function} is reduced to the symmetric case of the sum \eqref{bvp6}. To study more general integrals we will need
\begin{Prop}\label{self_reciprocal}
The functions of the form $f(x)=S\big(\tfrac{x}{\sqrt{2\pi}}\big)\cos\tfrac{x^2}{2}$, $g(x)=S\big(\tfrac{x}{\sqrt{2\pi}}\big)\sin\tfrac{x^2}{2}$, where
\[
S(-x)=S(x),\quad S(1-x)=-S(x),
\]
(i.e., $S(x)$ is a $2$-periodic function, even with respect to integer points, and odd with respect to half integer points) are self-reciprocal under cosine Fourier transform with eigenvalues $1$ and $-1$, respectively:
\[
F_c(x)=f(x),\quad G_c(x)=-g(x).
\]
\end{Prop}
\noindent{}\textit{Proof.} $S(x)$ has Fourier series representation
\[
S(x)=\sum_{n=0}^\infty c_n\cos\left(\pi x(2n+1)\right).
\]
Defining $h(x)=S\big(\tfrac{x}{\sqrt{2\pi}}\big)e^{ix^2/2}$, assuming that termwise integration of the series is possible, and calculating the integral using \eqref{gauss_integral}, we get
\begin{align*}
    H_c(x)=\sqrt{\frac{2}{\pi}}\sum_{n=0}^\infty c_n\cdot \int_0^\infty e^{{it^2}/{2}}\cos\Big(\sqrt{\tfrac{\pi}{2}}\,(2n+1)t\Big)\cos (xt) dt\\
=\sqrt{\frac{2}{\pi}}\sum_{n=0}^\infty c_n\cdot \sqrt{\frac{\pi i}{2}}\,e^{{-ix^2}/{2}-\pi i(2n+1)^2/4}\cos\Big(\sqrt{\tfrac{\pi}{2}}\,(2n+1)x\Big).
\end{align*}
To simplify the summand we note that $e^{-\pi i(2n+1)^2/4}=e^{-\pi i/4}$ for integer $n$, and get
\[
H_c(x)=e^{-ix^2/2}\sum_{n=0}^\infty c_n\cos\Big(\sqrt{\tfrac{\pi}{2}}\,(2n+1)x\Big)=e^{-ix^2/2}S\bigl(\tfrac{x}{\sqrt{2\pi}}\bigr),
\]
from which the claim follows. \qed 

\noindent{}\textit{Example.} Taking $S(x)=(-1)^{\lfloor x+\frac{1}{2}\rfloor}$ in Proposition \ref{self_reciprocal} gives two self-reciprocal functions
\[
(-1)^{\lfloor {x}/{\sqrt{2\pi}}+{1}/{2}\rfloor}\cos\tfrac{x^2}{2},\qquad (-1)^{\lfloor {x}/{\sqrt{2\pi}}+{1}/{2}\rfloor}\sin\tfrac{x^2}{2}.
\]
Combining these with generalized Plancherel's formula (\cite{titchmarsh}, 2.1.22)
\[
\sqrt{\alpha}\int_0^\infty f(x)g(\alpha x) dx=\sqrt{\beta}\int_0^\infty F_c(x)G_c(x\beta) dx,\quad \alpha\beta=1,
\]
and the self-reciprocal function \eqref{self_reciprocal_function1} leads to transformation formulas
\begin{equation}\label{plancherel_even}
    \sqrt{\alpha}\int_0^\infty\frac{\cos(\pi x^2)}{\cosh(\alpha x)}(-1)^{\lfloor x+\frac{1}{2}\rfloor}dx=\sqrt{\beta}\int_0^\infty\frac{\cos(\pi x^2)}{\cosh(\beta x)}(-1)^{\lfloor x+\frac{1}{2}\rfloor}dx,\quad \alpha\beta=\pi^2,
\end{equation}
\begin{equation}\label{plancherel_odd}
    \sqrt{\alpha}\int_0^\infty\frac{\sin(\pi x^2)}{\cosh(\alpha x)}(-1)^{\lfloor x+\frac{1}{2}\rfloor}dx+\sqrt{\beta}\int_0^\infty\frac{\sin(\pi x^2)}{\cosh(\beta x)}(-1)^{\lfloor x+\frac{1}{2}\rfloor}dx=0,\quad \alpha\beta=\pi^2.
\end{equation}
Symmetric case of \eqref{plancherel_odd} is the imaginary part of \eqref{integral_of_piecewise_function}. More generally: \textit{If} $F_c(x)=f(x)$, \textit{then}
\[
\int_0^\infty(-1)^{\lfloor x+\frac{1}{2}\rfloor}\,{\sin(\pi x^2)}f\big(\sqrt{2\pi}\,x\big)\,dx=0.
\]

\section{Trigonometric series for hypergeometric functions}

\subsection{}\label{fourier_series} Consider the Fourier series expansion
\begin{equation}\label{fourier}
    \frac{\pi e^{ia(\pi-x)}}{\sin(\pi a)}=\sum_{n\in\mathbb{Z}}\frac{e^{inx}}{n+a},\qquad 0<x<2\pi.
\end{equation}
Real and imaginary parts of \eqref{fourier} result in the formulas (\cite{gr1}, entries 1.445.5-6)
\begin{equation}\label{f1}
    \sum_{n=1}^\infty\left(\frac{\sin(n-a)x}{n-a}+\frac{\sin(n+a)x}{n+a}\right)=\pi-\frac{\sin(ax)}{a},\qquad 0<x<2\pi,
\end{equation}
\begin{equation}\label{f2}
    \sum_{n=1}^\infty\left(\frac{\cos(n+a)x}{n+a}-\frac{\cos(n-a)x}{n-a}\right)=\pi\cot(\pi a)-\frac{\cos(ax)}{a},\qquad |x|<2\pi.
\end{equation}
We used \eqref{f2} in Section \ref{sums_with_dirichlet_characters}, equation \eqref{piecewise_continuous_function}. Generalizations of \eqref{f1} and \eqref{f2} with higher powers of $(n\pm a)$ in the denominator are studied in (\cite{RNI}, Chapter 9, Entries 1 and 2). However, an interesting question is, what is the value of the truncated series in \eqref{fourier} where summation is over positive values of $n$ only. This question is answered by the less known formulas
\begin{equation}\label{3f2}
   \sum_{n=1}^\infty\left(\frac{\sin(n-a)x}{n-a}-\frac{\sin(n+a)x}{n+a}\right)={}_3F_2\left(\genfrac{}{}{0pt}{0}{{1}/{2},{1}/{2}+a,{1}/{2}-a}{{3}/{2},{3}/{2}};\sin^2\frac{x}{2}\right)\cdot 4a\sin\frac{x}{2},
\end{equation}
\begin{align}\label{4f3}
    \nonumber \sum_{n=1}^\infty\left(\frac{\cos(n+a)x}{n+a}+\frac{\cos(n-a)x}{n-a}\right)&={}_4F_3\left(\genfrac{}{}{0pt}{0}{1,1, 1+a,1-a}{2,2,{3}/{2} };\sin^2\frac{x}{2}\right)\cdot 2 a^2\sin^2\frac{x}{2}\\&-2\ln\left(2\sin\frac{x}{2}\right)+2\psi(1)-\psi(1-a)-\psi(1+a),
\end{align}
where $\psi$ denotes the digamma function. These formulas can be derived from Newton's formulas:
\begin{equation}\label{nf1}
    {}_2F_1\left(\genfrac{}{}{0pt}{0}{{1}/{2}+a,{1}/{2}-a}{{1}/{2}};\sin^2x\right)\cos x=\cos(2ax),
\end{equation}
\begin{equation}\label{nf2}
    {}_2F_1\left(\genfrac{}{}{0pt}{0}{{1}/{2}+a,{1}/{2}-a}{{3}/{2}};\sin^2x\right)\sin x=\frac{\sin(2ax)}{2a},
\end{equation}
\begin{equation}\label{nf3}
    {}_2F_1\left(\genfrac{}{}{0pt}{0}{a,-a}{{1}/{2}};\sin^2x\right)=\cos(2ax),
\end{equation}
(\cite{emot}, Entries 2.8(11), 2.8(12)). When $a\to 0$, formula \eqref{3f2} reduces to
\[
\sum_{n=0}^\infty\binom{2n}{n}\frac{\sin^{2n+1}x}{2^{2n}(2n+1)^2}=x\log|2\sin x|+\frac{1}{2}\sum_{n=1}^\infty\frac{\sin(2nx)}{n^2},
\]
given as Entry $16$ of Chapter 9 in \cite{RNI}. 

\subsection{} The formula
\begin{align}\label{gegenbauer}
    \nonumber\frac{4^{1-a-b}\sqrt{\pi}\,\Gamma(2a)\Gamma(2b)}{\Gamma(1-a-b)\Gamma(a+b+1/2)}&\,{}_2F_1\left(\genfrac{}{}{0pt}{0}{2a,2b}{a+b+1/2};\cos^2x\right)\\&=\sum_{n\in\mathbb{Z}}\frac{\Gamma(a+n/2)\Gamma(b+n/2)}{\Gamma(1-a+n/2)\Gamma(1-b+n/2)}\cos(2nx),\qquad x\in\mathbb{R}
\end{align}
was derived in \cite{cg} as a Fourier series expansion  of Gegenbauer functions (formula (1.2)). The formulas
\begin{equation}\label{12}
    \frac{\Gamma(1-a)\Gamma(1-b)}{2\Gamma(1-a-b)}\,{}_2F_1\left(\genfrac{}{}{0pt}{0}{a,b}{{1}/{2}};\cos^2x\right)=\frac{1}{2}+\sum_{n=1}^\infty\frac{(a)_n(b)_n}{(1-a)_n(1-b)_n}\cos(2nx),\qquad x\in\mathbb{R}
\end{equation}
\begin{equation}\label{32}
    \frac{\Gamma(2-a)\Gamma(2-b)}{2\Gamma(2-a-b)}\,{}_2F_1\left(\genfrac{}{}{0pt}{0}{a,b}{{3}/{2}};\cos^2x\right)\cos x=\sum_{n=0}^\infty\frac{(a)_n(b)_n}{(2-a)_n(2-b)_n}\cos(2n+1)x,\qquad x\in\mathbb{R}
\end{equation}
follow from \eqref{gegenbauer} and the quadratic transformation formulas 2.11(7) and 2.11(9) in \cite{emot}. A more direct method is to use the general approach of \cite{macrobert} (see also \cite{bajpai}).

Another trigonometric expansion formula is found in \cite{emot} (formula 3.5(2)), which we write here as
\begin{align}\label{trig_expansion}
    \nonumber\frac{\sqrt{\pi}\,\Gamma(1/2+a)}{2\Gamma(1+a-c)\Gamma(c)}&{}_2F_1\left(\genfrac{}{}{0pt}{0}{a,1-a}{c};\sin^2x\right)(2\sin x)^{2c-2}\\
    &=\sum_{k=0}^\infty\frac{(1+a-c)_k(3/2-c)_k}{k!(1/2+a)_k}\sin\big[2(2k+1+a-c)x\big],\qquad 0<x<\frac{\pi}{2}
\end{align}
The hypergeometric function appearing on the left hand side of \eqref{trig_expansion} is the so called Ferrers function up to an elementary function prefactor. More recently, \eqref{trig_expansion} was considered in \cite{volkmer}.
Note that the types of hypergeometric functions on the left hand sides of \eqref{gegenbauer} and \eqref{trig_expansion} are related by the quadratic transformation formula (\cite{emot}, Entry 2.11(2))
\[
{}_2F_1\left(\genfrac{}{}{0pt}{0}{2a,2b}{a+b+1/2};\sin^2x\right)={}_2F_1\left(\genfrac{}{}{0pt}{0}{a-b+1/2,b-a+1/2}{a+b+1/2};\sin^22x\right),\qquad |x|<\frac{\pi}{4}.
\]

Particular case of \eqref{trig_expansion}, namely $a=1/2$, $c=1$,
\begin{equation}\label{elliptic_fourier}
    K(\sin x)=\frac{\pi}{2}\cdot{}_2F_1\left(\genfrac{}{}{0pt}{0}{{1}/{2},{1}/{2}}{1};\sin^2x\right)=\pi\sum_{n=0}^\infty\frac{\left({1}/{2}\right)_n^2}{(n!)^2}\sin(4n+1)x,\qquad 0<x<\frac{\pi}{2},
\end{equation}
was derived by Tricomi \cite{tricomi1}, \cite{tricomi2} (also reproduced in \cite{emot}, Entry 13.8(8)), using a connection between Legendre expansions, Abel transform and triginometric Fourier series. More recently, the same connection appeared in the paper \cite{mv1} in the context of numerical algorithms, and was extended to Gegenbauer polynomials in \cite{mv2}. \eqref{elliptic_fourier} was used recently in \cite{bbbg} to calculate certain integrals of elliptic integrals.

\subsection{} The following formula which gives a Fourier series expansion of a generalized hypergeometric function ${}_3F_2$ of the argument $\sin^2x$ is motivated by the formulas above:
\begin{Prop} If $\,\operatorname{Re} (a+b)<1$, and $b$ is not a positive integer, then
    \begin{align}\label{hypergeometric1}
    \nonumber &{}_3F_2\left(\genfrac{}{}{0pt}{0}{1,a,b}{(a+b)/2,{(1+a+b)}/{2}};\sin^2x\right)\cos x=\frac{a+b-1}{b-1}\sum_{n=0}^\infty \frac{(a)_n}{(2-b)_n}\cos(2n+1)x\\
    &\phantom{..............................................}+\frac{\Gamma(a+b)\Gamma(1-b)}{\Gamma(a)(2\sin x)^{a+b-1}}\sin\left(\tfrac{\pi}{2}(a+b)+x(b-a)\right),\qquad 0<x<\frac{\pi}{2}.
\end{align}
\end{Prop}

\noindent{}\textit{Proof.} Starting from \eqref{nf1} written in the form
\[
\sum_{k=0}^\infty\frac{(1/2+s)_k(1/2-s)_k}{k!(1/2)_k}\sin^{2k}x=\frac{\cos(2sx)}{\cos x},
\]
we multiply it by
\begin{equation}\label{gamma}
    \Gamma(1/2-s)\Gamma(1/2+s)\Gamma(c-s)\Gamma(d+s)
\end{equation}
to obtain
\[
\sum_{k=0}^\infty\Gamma(1/2+k+s)\Gamma(1/2+k-s)\Gamma(c-s)\Gamma(d+s)\frac{\sin^{2k}x}{k!(1/2)_k}=\frac{\pi\Gamma(c-s)\Gamma(d+s)\cos(2sx)}{\cos x\cos(\pi s)},
\]
and then integrate over $s$ along the contour $C$ parallel to the imaginary axes that separates the increasing set of poles of \eqref{gamma} from the decreasing set of poles. The integral on the left hand side can be calculated using Barnes' first lemma \cite{aar}. The integral on the right hand side is convergent because of the condition $0<x<{\pi}/{2}$. It can be calculated by closing the contour $C$ on the left half plane, picking up residues at $-1/2-n$ and $-n-d$ with $n\ge 0$, as described in \cite{aar}. The sum over residues at $-1/2-n$ results in the first term on the right hand side of \eqref{hypergeometric1}. It turns out, that the sum over the residues $-n-d$ is summable by the binomial theorem, and it results in the second term on the right hand side of \eqref{hypergeometric1}.\qed

\eqref{hypergeometric1} is essentially a three-term quadratic transformation formula for ${}_3F_2$ when arguments of the two of the functions ${}_3F_2$ lie on the unit circle. Formulas \eqref{12}, \eqref{32}, \eqref{hypergeometric1} are generalizations of \eqref{nf1}, \eqref{nf2}, \eqref{nf3} in the following sense. Setting $b=1-a$ and changing $a$ to $1/2-a$ in \eqref{hypergeometric1} immediately yields \eqref{nf1}. After setting $b=-a$ in \eqref{12}, one may notice that the series on right hand side is Fourier series expansion of $\frac{\pi a}{2\sin(\pi a)}\cos (\pi-2x)a$ (section \ref{fourier_series}). Similarly, \eqref{nf2} follows from \eqref{32} when $b=1-a$ and redefinition of parameters.

The series on the right hand side of \eqref{hypergeometric1} is a truncation of the bilateral series
\begin{equation}\label{bilateral}
    \sum_{n\in\mathbb{Z}} \frac{(a)_n}{(2-b)_n}\cos(2n+1)x.
\end{equation}
Due to Riemann's form of the binomial theorem (when the argument lies on the unit circle), the series \eqref{bilateral} has a closed form in terms of trigonometric functions (e.g., \cite{osler}, and references therein). While the truncated series does not have a simple closed form, \eqref{hypergeometric1} shows that it is related to a hypergeometric series with argument $\sin^2x$, similar to the sitution encountered in subsection \ref{fourier_series}.

\section{Inverse tangent integral and associated infinite products}

\subsection{}  Consider the problem submitted by Ramanujan to the Journal of the Indian Mathematical Society and solved by him in the article \cite{ramanujan} (also \cite{bck}, page 31; \cite{RNV}, chapter 37, entry 30): {\it Let $\alpha>0$ and $0<\beta<1$, with}
\[
\frac{1}{2}\pi\alpha=\log\tan\left\{\frac{1}{4}\pi(1+\beta)\right\},
\]
\textit{then}
\begin{equation}\label{inf_product}
    \left(\frac{1^2+\alpha^2}{1^2-\beta^2}\right)\left(\frac{3^2-\beta^2}{3^2+\alpha^2}\right)^3\left(\frac{5^2+\alpha^2}{5^2-\beta^2}\right)^5\ldots =e^{\pi\alpha\beta/2}.
\end{equation}
Here, the condition on $\alpha$ and $\beta$ can be cast in the symmetrical form
\[
\tanh\Big(\frac{\pi\alpha}{4}\Big)=\tan\Big(\frac{\pi\beta}{4}\Big),
\]
or equivalently
\[
\cosh\Big(\frac{\pi\alpha}{2}\Big)\cos\Big(\frac{\pi\beta}{2}\Big)=1,
\]
and the product can be written succinctly as
\[
\prod_{n=1}^\infty  \left(\frac{n^2+\alpha^2}{n^2-\beta^2}\right)^{n\chi_4(n)}.
\]

The logarithm of the infinite product in \eqref{inf_product} can be expressed through the inverse tangent integral
\[
\operatorname{Ti}_2(z)=\int_0^{z}\frac{\tan^{-1} t}{t}\,dt,
\]
as was shown in \cite{ramanujan}. There is a trigonometric series expansion for the inverse tangent integral (\cite{ramanujan}; \cite{RNI}, chapter 9, Entry 17; also mentioned in \cite{lewin}, section 2.4.2):
\[
\operatorname{Ti}_2(\tan x)=\sum_{n=0}^\infty \frac{(-1)^n\tan^{2n+1}x}{(2n+1)^2}=x\log|\tan x|+\frac{1}{2}\sum_{n=0}^\infty\frac{\sin(4n+2)x}{(2n+1)^2},\qquad |x|<\frac{\pi}{4}.
\]
Thus, the logarithm of the infinite product in \eqref{inf_product} can be expressed in terms of dilogarithm function. More generally:
\begin{Prop} Let $\chi_k(n)$ be an odd Dirichlet character modulo $k$ (or any $k$-periodic odd sequence \eqref{odd_periodic_sequence}) and $0<\beta<1$. Then
    \begin{equation}\label{chi_product}
    \prod_{n=1}^\infty\left(1-\frac{\beta^2}{n^2}\right)^{n\chi_k(n)}=\exp\left(-\frac{\pi}{k}\sum_{j=1}^{k-1}\chi_k(j)\int_0^\beta x\cot\frac{\pi(x+j)}{k}\,dx\right).
\end{equation}
\end{Prop}
\noindent{}Thus, the infinite product on the left hand side of \eqref{chi_product} can be expressed in terms of exponential of a linear combination of dilogarithms (after integration by parts, the integrals in \eqref{chi_product} are expressed in terms of Clausen function, which in turn can be expressed in terms of dilogarithms). Closely related to \eqref{chi_product} is the integral representation for the logarithm of the ratio of Barnes $G$-functions
\[
\log\frac{G(1+z)}{G(1-z)}=z\log(2\pi)-\int_0^z\pi x\cot(\pi x)\, dx
\]
(\cite{ww}, chapter 12, example 49; also, equation (33) in \cite{adamchik}; or the recent historical account \cite{neretin}). 

\subsection{}
The considerations below were motivated by a search for a finite version of \eqref{inf_product}.

\begin{Lemma}\label{lemma} Let $m$ be a non-negative integer, $\alpha>0$, $0<\beta<\pi/2$. Then
\begin{align}\label{auxilary_formula}
\nonumber&\int_0^\alpha \frac{\sinh\frac{x}{2 m+1}}{\cosh x} \, dx+\int_0^\beta \frac{\sin\frac{y}{2 m+1}}{\cos y} \, dy\\&=(-1)^m \ln \frac{\cosh\frac{\alpha}{2 m+1}}{\cos \frac{\beta}{2 m+1}}+\sum _{k=0}^{m-1} (-1)^k \sin \frac{\pi  (2 k+1)}{2 (2 m+1)}\cdot\log\frac{\sin ^2\frac{\pi  (2 k+1)}{2 (2 m+1)}+\sinh ^2\frac{\alpha}{2 m+1}}{\sin ^2\frac{\pi  (2 k+1)}{2 (2 m+1)}-\sin ^2\frac{\beta}{2 m+1}}.
\end{align}
\end{Lemma}

\noindent {\it Proof.}  It is quite straightforward to verify the partial fractions expansion
\begin{equation}\label{part_frac} 
\frac{2m+1}{\cosh \left((2m+1) x\right) }=\sum _{|k|\le m} (-1)^{m-k}\frac{\cos\frac{\pi k}{2m+1}\cosh x}{\sinh^2x+ \cos^2\frac{\pi  k}{2m+1}}.
\end{equation}
By multiplying this formula by $\sinh x$ and integrating with respect to $x$ from $0$ to $\alpha/(2m+1)$ one obtains evaluation of the first integral in \eqref{auxilary_formula}, that is \eqref{auxilary_formula} with $\beta=0$. Similarly, considering the companion partial fraction expansion
\[
\frac{2m+1}{\cos \left((2m+1) y\right)}=\sum _{|k|\le m} (-1)^{m-k}\frac{\cos\frac{\pi k}{2m+1}\cos y}{\cos^2\frac{\pi  k}{2m+1}-\cos^2y},
\]
one evaluates the second integral in \eqref{auxilary_formula} (the condition $\beta<\pi/2$ ensures that $\cos^2\frac{\pi  k}{2m+1}-\cos^2y>0$ for all $|k|\le m$). Then, their sum is the quantity on the right hand side of \eqref{auxilary_formula}. \qed

One can deduce \eqref{inf_product} from the Lemma \ref{lemma} as follows. Assume that $\alpha$ and $\beta$ in Lemma \ref{lemma} are related by $\cosh(\alpha)\cos(\beta)=1$. We multiply \eqref{auxilary_formula} by $(2m+1)$ and then take the limit $m\to\infty$. It may be easily verified that when $\alpha$ and $\beta$ are subject to the constraint $\cosh(\alpha)\cos(\beta)=1$, then
\[
\frac{d\alpha}{\cosh\alpha}=d\beta,\quad \frac{d\beta}{\cos\beta}=d\alpha.
\]
As a result
\[
\int_0^\alpha \frac{x dx}{\cosh x}+\int_0^\beta \frac{ydy}{\cos y}=\alpha\beta,
\]
when $\cosh(\alpha)\cos(\beta)=1$. Hence, the left hand side of \eqref{auxilary_formula} multiplied by $(2m+1)$ tends to $\alpha\beta$ in the limit $m\to\infty$. The right hand side multiplied by $(2m+1)$ tends to a logarithm of a certain infinite product, which after rescaling $\alpha$ and $\beta$ by a factor of $\pi/2$ coincides with the infinite product in \eqref{inf_product}. Exponentiation of both sides then completes the proof of \eqref{inf_product}. Thus, we were able to give an elementary proof of \eqref{inf_product}.

\subsection{} Is it possible that an identity of the similar type as \eqref{inf_product} exists for Dirichlet characters other than $\chi_4$? Consider the elementary integral
\[
\int_0^\alpha\frac{\sinh (x)}{\sinh (3x)}\,dx=\frac{1}{\sqrt{3}}\arctan \left(\frac{\tanh\alpha}{\sqrt{3}}\right).
\]
It follows from this that when $\frac{\pi}{3}<\beta<\frac{\pi}{2}$ and $\tanh(\alpha)\tan(\beta)=\sqrt{3}$
\[
\int_0^\alpha\frac{x\sinh (x)}{\sinh (3x)}\,dx+\int_\beta^{\pi/2}\frac{y\sin (y)}{\sin (3y)}\,dy=-\frac{\alpha\beta}{\sqrt{3}}.
\]
Alternatively, one can check this by differentiation. Moreover
\[
\sum_{n=1}^\infty \left(\frac{n}{3}\right)\frac{n}{x^2+n^2} =\frac{\pi}{\sqrt{3}}\frac{\sinh\left(\frac{\pi x}{3}\right)}{\sinh (\pi x)}, 
\]
where $\left(\tfrac{n}{3}\right)$ is the Legendre symbol mod $3$. Thus, one arrives at the following result:
\begin{Prop}\label{inf_prod1}
    If $\alpha>0$, $1<\beta<3/2$ and $\tanh\big(\frac{\pi\alpha}{3}\big)\tan\big(\frac{\pi\beta}{3}\big)=\sqrt{3}$, then
\[
\prod_{n=1}^\infty\left(\frac{(n^2+\alpha^2)(n^2-\beta^2)}{n^2(n^2-9/4)}\right)^{n\left(\tfrac{n}{3}\right)}=e^{-2\pi\alpha\beta/3}.
\]
\end{Prop}
\subsection{} Similarly, when $0<\beta<\frac{\pi}{4}$ and $\tanh(\alpha)=\sqrt{3}\,\tan(\beta)$
\[
\int_0^\alpha\frac{x\sinh (x)}{\sinh (3x)}\,dx+\int_0^\beta\frac{y\sin (y)}{\sin (3y)}\,dy+\frac{1}{2}\int_0^{2\beta}\frac{y\sin (y)}{\sin (3y)}\,dy=\frac{\alpha\beta}{\sqrt{3}}.
\]
This results in two formulas that are dual to each other:
\begin{Prop}\label{inf_prod2} Let $\alpha>0$.
\begin{enumerate}[label=\textnormal{(\roman*)}, wide=0pt]
\item If $0<\beta<1/2$ and $\tanh\big(\frac{\pi\alpha}{3}\big)=\sqrt{3}\,\tan\big(\frac{\pi\beta}{3}\big)$, then
    \[
\prod_{n=1}^\infty\left(\frac{n^2(n^2+\alpha^2)^2}{(n^2-\beta^2)^2(n^2-4\beta^2)}\right)^{n\left(\tfrac{n}{3}\right)}=e^{4\pi\alpha\beta/3}.
\]
\item If $0<\beta<1$ and $\tan\big(\frac{\pi\beta}{3}\big)=\sqrt{3}\,\tanh\big(\frac{\pi\alpha}{3}\big)$, then
\[
\prod_{n=1}^\infty\left(\frac{n^2(n^2-\beta^2)^2}{(n^2+\alpha^2)^2(n^2+4\alpha^2)}\right)^{n\left(\tfrac{n}{3}\right)}=e^{-4\pi\alpha\beta/3}.
\]
\end{enumerate}
\end{Prop}

\subsection{} When $\frac{\pi}{6}<\beta<\frac{\pi}{2}$ and $\tanh(\alpha)\tan(\beta)=1/\sqrt{3}$
\[
\int_0^\alpha\frac{x\cosh (x)}{\cosh (3x)}\,dx+\int_\beta^{\pi/2}\frac{y\cos (y)}{\cos (3y)}\,dy=-\frac{\alpha\beta}{\sqrt{3}}.
\]
Also
\[
\sum_{n=1}^\infty \chi_{6}(n)\frac{n}{x^2+n^2} =\frac{\pi}{6\sqrt{3}}\frac{\cosh\left(\frac{\pi x}{6}\right)}{\cosh \left(\frac{\pi x}{2}\right)},
\]
where
\begin{equation}
    \chi_{6}(n)=\begin{cases}
         \phantom{-}1, & n\equiv 1~~ (\operatorname{mod} 6),\\
         -1, & n\equiv 5~~ (\operatorname{mod} 6), \\
         \phantom{-} 0, & \mathrm{otherwise},
    \end{cases}
\end{equation}
is the odd Dirichlet character modulo $6$. Hence
\begin{Prop}\label{inf_prod3} If $\alpha>0$, $1<\beta<3$ and $\tanh\big(\frac{\pi\alpha}{6}\big)\tan\big(\frac{\pi\beta}{6}\big)=1/\sqrt{3}$, then 
\[
\prod_{\substack{n=1\\ n\neq 3}}^\infty\left(\frac{(n^2+\alpha^2)(n^2-\beta^2)}{n^2(n^2-9)}\right)^{n\chi_{6}(n)}=e^{-\pi\alpha\beta/3}.
\]
\end{Prop}

\subsection{} Due to the formula \eqref{chi_product}, Propositions \ref{inf_prod1} - \ref{inf_prod3} can also be proved directly using functional equations for the dilogarithm function, though calculations in this case are more cumbersome. Also, previous results can be generalized to include an additional continuous parameter.
\begin{Prop} Let $\theta\in(0,\pi)$ and $\alpha>0$.
\begin{enumerate}[label=\textnormal{(\roman*)}, wide=0pt]
\item If $0<\beta<\frac{\theta}{2}$ and $\tanh(\alpha)=\cot\left(\frac{\theta}{2}\right)\tan(\beta)$, then
\[
\int_0^\alpha \frac{x dx}{\cosh(2x)+\cos\theta}+\int_0^\beta \frac{ydy}{\cos(2y)-\cos\theta}=\frac{\alpha\beta}{\sin\theta}.
\]
\item If $\frac{\pi-\theta}{2}<\beta<\frac{\pi}{2}$ and 
$\tanh(\alpha)\tan(\beta)=\cot\left(\frac{\theta}{2}\right)$, then
\[
\int_0^\alpha \frac{x dx}{\cosh(2x)+\cos\theta}+\int_\beta^{\pi/2} \frac{ydy}{\cos(2y)+\cos\theta}=-\frac{\alpha\beta}{\sin\theta}.
\]
\end{enumerate}
\end{Prop}
\noindent{}To convert this into product form, one can use the partial fractions expansion
\[
\frac{1}{\cosh(\pi x)+\cos(\pi\theta)}=\frac{2}{\pi\sin(\pi\theta)}\sum_{n=0}^\infty\left(\frac{2 n+1-\theta}{x^2+(2 n+1-\theta)^2}-\frac{2 n+1+\theta}{x^2+(2 n+1+\theta)^2}\right).
\]

\subsection{} Let $\sin_px$ and $\cos_px$ denote the generalized trigonometric functions
defined as the solution to the differential equation \cite{yin}
\[
\frac{d}{dx}\sin_px=\cos_px,\qquad \sin_p^px+\cos_p^px=1,\qquad \sin_p0=0,
\]
when $x\in[0;\pi_p/2]$, where
\[
\pi_p=2\int_0^1\frac{dt}{(1-t^p)^{1/p}}.
\]
The generalized hyperbolic functions are defined as a solution to the differential equation \cite{yin}
\[
\frac{d}{dx}\sinh_px=\cosh_px,\qquad \cosh_p^px-\sinh_p^px=1,\qquad \sinh_p0=0.
\]
There is a certain duality between the generalized trigonometric and hyperbolic functions. Such a duality has been explicitly stated in \cite{mt}. On the other hand, the generalized Gudermannian function was defined by Neuman \cite{neuman} as
\[
\operatorname{gd}_p(x)=\int_0^x \frac{dt}{\cosh_p^{p-1} t},
\]
where he has shown that the inverse function of the generalized Gudermannian function is
\[
\operatorname{gd}^{-1}_p(x)=\int_0^x \frac{dt}{\cos_p^{p-1} t}.
\]
If $\beta=\operatorname{gd}_p(\alpha)$ (or equivalently $\alpha=\operatorname{gd}^{-1}_p(\beta)$) it follows from $\operatorname{gd}'_p(\alpha)\cdot \left[\operatorname{gd}^{-1}_p(\beta)\right]'=1$ and Neuman's result, that $\alpha$ and $\beta$ are subject to the condition $\cosh_p\alpha\cdot \cos_p\beta=1$, and vice versa. This implies: 
\begin{Prop}
If $\,\alpha>0$, $0<\beta<\pi_p/2$ and $\cosh_p\alpha\cdot \cos_p\beta=1$, then
\begin{equation}\label{gtf}
    \int_0^\alpha \frac{xdx}{\cosh_p^{p-1} x}+\int_0^\beta \frac{ydy}{\cos_p^{p-1} y}=\alpha\beta.
\end{equation}  
\end{Prop}

In general, there are no simple partial fractions expansion formulas for generalized trigonometric functions. For example, it has been shown in \cite{egl} that for certain set of parameters generalized trigonometric functions are expressed in terms of Jacobi elliptic functions, and thus they are doubly periodic. This means their partial fractions expansion is a double series, similar to Weierstrass elliptic functions. This probably means that there are no formulas similar to \eqref{inf_product} resulting from \eqref{gtf} unless $p=2$.

\end{document}